\documentclass{amsart}

\usepackage{amsfonts,amssymb,amscd,amsmath,enumerate,verbatim}

\usepackage[all]{xy}
\SelectTips{eu}{} 
\xyoption{curve}

\newcommand{\lra}{\longrightarrow}

\newcommand{\wt}{\widetilde}

\newcommand{\col}{\colon}

\newcommand\dd{\partial}
\newcommand{\hh}[1]{\operatorname{H}(#1)}
\newcommand{\HH}[2]{\operatorname{H}_{#1}(#2)}
\newcommand{\Tor}[4]{\operatorname{Tor}_{#1}^{#2}(#3,#4)}

\newcommand\Ker{\operatorname{Ker}}

\newcommand{\D}{\operatorname{D}}
\newcommand{\Tom}{\operatorname{Hom}}

\newcommand{\cond}[3][R]{\left(#3:_{#1}#2\right)}
\newcommand\depth{\operatorname{depth}}
\newcommand\height{\operatorname{height}}
\newcommand\length{\operatorname{length}}
\newcommand{\rank}{\operatorname{rank}}

\newcommand\fm{{\mathfrak m}}
\newcommand\fp{{\mathfrak p}}
\newcommand\bsx{{\boldsymbol x}}
\newcommand\BN{{\mathbb N}}

\theoremstyle{plain}

\newtheorem{theorem}{Theorem}[section]

\newtheorem{proposition}[theorem]{Proposition}
\newtheorem{lemma}[theorem]{Lemma}

\newtheorem{corollary}[theorem]{Corollary}

\theoremstyle{definition}

\newtheorem{remark}[theorem]{Remark}

\newtheorem{chunk}[theorem]{}

\theoremstyle{remark}

\numberwithin{equation}{theorem}

\begin{document}

\title[Regular local rings]
{A criterion for regularity of local rings}

\date{\today}

\author[T.~Bridgeland]{Tom Bridgeland} 
\address{Department of Pure Mathematics, University of Sheffield, Sheffield S3 7RH, U.K.}
\email{t.bridgeland@shef.ac.uk}

\author[S.~Iyengar]{Srikanth Iyengar} 
\address{Department of Mathematics, University of Nebraska, Lincoln, NE 68588, U.S.A.}
\email{iyengar@math.unl.edu}

\thanks {S.I. was partly supported by NSF grant DMS 0442242}

\begin{abstract}
  It is proved that a noetherian commutative local ring $A$ containing a field is regular
  if there is a complex $M$ of free $A$-modules with the following properties: $M_i=0$ for
  $i\notin[0,\dim A]$; the homology of $M$ has finite length; $\HH 0M$ contains the
  residue field of $A$ as a direct summand.  This result is an essential component in the
  proofs of the McKay correspondence in dimension 3 and of the statement that threefold
  flops induce equivalences of derived categories.
\end{abstract}

\maketitle

\section{Introduction}
\label{Introduction}

Let $(A,\fm,k)$ be a local ring; thus $A$ is a noetherian commutative local ring, with
maximal ideal $\fm$ and residue field $k$.  Let $M\col 0\to M_d\to \cdots\to M_0\to 0$ be
a complex of free $A$-modules with $\length_A\hh M$ finite and non-zero.  The New
Intersection Theorem \cite{Pr} yields $d\geq \dim A$.  

In this article we prove the following result, which is akin to Serre's theorem that a
local ring is regular when its residue field has a finite free resolution.

\begin{theorem}
  Assume $A$ contains a field or $\dim A\leq 3$.  If $d=\dim A$ and $k$ is a direct
  summand of $\HH 0M$, then $\HH iM=0$ for $i\ge 1$, and the ring $A$ is regular.
\end{theorem}

This result is contained in Theorem \ref{really main}. A restricted version of such a
statement occurs as \cite[(4.3)]{BM}; however, as is explained in Remark \ref{compare},
the proof of \emph{loc.~cit.} is incorrect.

In the remainder of the introduction, the field $k$ is algebraically closed, schemes over
$k$ are of finite type, and the points considered are closed points. We write $\D(Y)$ for
the bounded derived category of coherent sheaves on a scheme $Y$, and $k_y$ for the
structure sheaf of a point $y\in Y$.  Given the theorem above, arguing as in
\cite[\S5]{BM}, one obtains the following corollary.

\begin{corollary}
  Let $Y$ be an irreducible scheme of dimension $d$ over $k$, and let $E$ be an object of
  $\D(Y)$. Suppose there is a point $y_0\in Y$ such that $k_{y_0}$ is a direct summand of
$\HH 0E$ and
\[
\Tom_{\D(Y)}(E,{k_y[i]}) = 0 \text{ unless }y=y_0\text{ and }0\leq i\leq d.
\]
Then $Y$ is non-singular at $y_0$ and $E\cong \HH 0E$ in $\D(Y)$. \qed
\end{corollary}

This result enables one to show that certain moduli spaces are non-singular and give rise
to derived equivalences; see \cite[(6.1)]{BM}. This has proved particularly effective in
dimension three, and is an essential component in the proofs in \cite{B,BKR}.

\section{Proof of the main theorem}
\label{Proof}
This section is dedicated to a proof Theorem \ref{really main}. The book of Bruns and
Herzog \cite{BH} is our standard reference for the notions that appear here.

Let $(A,\fm,k)$ be a local ring, and let $C$ be an $A$-module; it need not be finitely
generated.  A sequence $\bsx = x_1,\dots,x_n$ is \emph{$C$-regular} if $\bsx C \ne C$ and
$x_{i}$ is a non-zero divisor on $C/(x_1,\dots,x_{i-1})C$ for each $1\leq i\leq n$.  An
$A$-module $C$ is \emph{big Cohen-Macaulay} if there is a system of parameters for $A$
that is $C$-regular.  If every system of parameters of $A$ is $C$-regular, then $C$ is
said to be \emph{balanced}. Any ring that has a big Cohen-Macaulay module has also one
that is balanced; see \cite[(8.5.3)]{BH}.  Big Cohen-Macaulay modules were introduced by
Hochster \cite{Ho}, who constructed them when $A$ contains a field, and, following the
recent work of Heitmann, also when $\dim A\leq 3$; see \cite{Ho1}.

The result below is contained in the proof of \cite[(1.13)]{EG} by Evans and Griffith, see
also \cite[(3.1)]{mz}, so only an outline of an argument is provided; it follows
the discussion around \cite[(9.1.7)]{BH}; see also \cite[(3.4)]{mz}.

\begin{lemma}
\label{depth sensitivity}
Let $M$ be a complex of free $A$-modules with $M_i=0$ for $i\notin[0,\dim A]$, the
$A$-module $\HH 0M$ finitely generated, and $\length_A\HH iM$ finite for each $i\geq
1$. 

If $C$ is a balanced big Cohen-Macaulay module, then $\HH i{M\otimes_AC}=0$ for $i\ge
1$.
\end{lemma}

\begin{proof}[Sketch of a proof]
Replacing $M$ by a quasi-isomorphic complex we may assume each $M_i$ is finite free.
Choose a basis for each $M_i$ and let $\phi_i$ be the matrix representing the differential
$\dd_i \col M_i\to M_{i-1}$. Set $r_i = \sum_{j=i}^n(-1)^{j-i}\rank_A M_j$, and let
$I_{r_i}(\phi_i)$ be the ideal generated by the $r_i\times r_i$ minors of $\phi_i$.  Fix
an integer $1\leq i\leq d$, and let $\fp$ be a prime ideal of $A$ with $\height
I_{r_i}(\phi_i)=\height \fp$. If $\height \fp=\dim A$, then $\height I_{r_i}(\phi_i)=\dim
A \geq d\geq i$, where the first inequality holds by hypotheses.  If $\height \fp < \dim A$,
then $\HH j{M\otimes_AA_\fp}=0$ for $j\geq 1$, by hypotheses. Therefore
\[
\height I_{r_i}(\phi_i) = \height(I_{r_i}(\phi_i)\otimes_A A_\fp) = 
\height I_{r_i}(\phi_i\otimes_A A_\fp) \geq i 
\]
where the inequality is comes from the Buchsbaum-Eisenbud acyclicity criterion
\cite[(9.1.6)]{BH}.  Thus, no matter what $\height\fp$ is, one has
\[
\dim A - \dim A/I_{r_i}(\phi_i) \geq \height I_{r_i}(\phi_i) \geq i
\]
Thus, $I_{r_i}(\phi_i)$ contains a sequence $\bsx=x_1,\dots,x_i$ that extends to a full
system of parameters for $R$. Since $C$ is balanced big Cohen-Macaulay module, the
sequence $\bsx$ is $C$-regular, so another application of \cite[(9.1.6)]{BH} yields the
desired result.
\end{proof}

The following elementary remark is invoked twice in the arguments below.

\begin{lemma}
\label{tor vanishing}
Let $R$ be a commutative ring and $U,V$ complexes of $R$-modules with $U_i = 0 = V_i$ for
each $i<0$. If each $R$-module $V_i$ is flat and $\HH 1{U\otimes_RV}=0$, then $\HH 1{{\HH
0U}\otimes_RV}=0$.
\end{lemma}

\begin{proof}
Let $\wt U=\Ker(U\to \HH 0U)$; evidently, $\HH i{\wt U}=0$ for $i<1$. Since $-\otimes_RV$
preserves quasi-isomorphisms, $\HH i{\wt U\otimes_RV}=0$ for $i<1$.  The long exact
sequence that results from the short exact sequence of complexes $0\to \wt U\to U\to \HH
0U\to 0$ thus yields a surjective homomorphism
\[
\HH 1{U\otimes_RV} \lra \HH 1{\HH 0U\otimes_RV} \lra 0
\]  
This justifies the claim.
\end{proof}

The proposition below is immediate when $C$ is finitely generated: $\Tor
1ACk=0$ implies  $C$ is free.  It contains \cite[(2.5)]{Sc}, due to
Schoutens, which deals with the case when $C$ is a big Cohen-Macaulay algebra.

\begin{proposition}
\label{cm:criteria}
Let $(A,\fm,k)$ be a local ring and $C$ an $A$-module with $\fm C\ne C$.  If $\Tor
1ACk=0$, then each $C$-regular sequence is $A$-regular.
\end{proposition}

\begin{proof}
First we establish that for any ideal $I$ in $A$ one has $\cond C{IC}=I$.

Indeed, consider first the case where the ideal $I$ is $\fm$-primary ideal.
  
Let $F$ be a flat resolution of $C$ as an $A$-module, and let $V$ be a flat resolution of
$k$, viewed as an $A/I$-module. Therefore $\HH 1{F\otimes_AV}=\Tor 1ACk=0$.  The $A$
action on $k$ factors through $A/I$, so $F\otimes_A V\cong
(F\otimes_AA/I)\otimes_{A/I}V$. Thus, applying Lemma \ref{tor vanishing} with
$U=C\otimes_AA/I$ one obtains
\[
\Tor 1{A/I}{C/IC}k = \HH 1{\HH 0 U\otimes_{A/I}V}=0
\]
The ring $A/I$ is artinian and local, with residue field $k$, thus $\Tor 1{A/I}{C/IC}k=0$
implies that the $A/I$-module $C/IC$ is free; see, for instance, \cite[(22.3)]{Ma}. Moreover,
$C/IC$ is non-zero as $\fm C\ne C$.  Thus $\cond C{IC}=I$, as desired.

For an arbitrary ideal $I$, evidently $I\subseteq \cond C{IC}$. The reverse inclusion
follows from the chain:
\begin{align*}
\cond C{IC} \subseteq \bigcap_{n\in\BN}\cond C{(I+\fm^n)C} 
            =  \bigcap_{n\in\BN} (I+\fm^n)
            = I
\end{align*}
where the first equality holds because each ideal $(I+\fm^n)$ is a $\fm$-primary, while
the second one is by the Krull Intersection Theorem. This settles the claim.

Let $x_1,\dots,x_m$ be a regular sequence on $C$.  Fix an integer $1\leq i\leq m$, and set
$I=(x_1,\dots,x_{i-1})$.  For any element $r$ in $A$, the first and the second
implications below are obvious:
\begin{align*}
rx_{i}  \in I  \implies rx_{i}C\subseteq IC 
             \implies x_{i}(rC)\subseteq IC 
             \implies rC\subseteq IC 
             \implies r\in I
\end{align*}
The third implication holds because $x_i$ is regular on $C/IC$, and the last one is by the
claim established above. Thus, $x_{i}$ is a non-zero divisor on $A/I$, that is to say, on
$A/(x_1,\dots,x_{i-1})$.  Since this holds for each $i$, we deduce that the sequence
$\bsx$ is regular on $A$, as desired.
\end{proof}

The result below contains the theorem stated in the introduction.

\begin{theorem}
\label{really main}
Let $(A,\fm,k)$ be a local ring, $M$ a complex of free $A$-modules with $M_i=0$ for
$i\notin[0,\dim A]$, the $A$-module $\HH 0M$ finitely generated, and $\length_A\HH iM$
finite for $i\geq 1$.  Assume $A$ has a big Cohen-Macaulay module.
If $k$ is a direct summand of $\HH 0M$, then the local ring $A$ is regular.
\end{theorem}

\begin{proof}
  Let $C$ be a big Cohen-Macaulay $A$-module; we may assume that $C$ is balanced. Let
  $d=\dim A$, and let $\bsx=x_1,\dots,x_d$ be a system of parameters for $A$ that is a
  regular sequence on $C$. In particular, $\bsx C\ne C$, and hence $\fm C\ne C$, since
  $\bsx$ is $\fm$-primary.

Let $V$ be a flat resolution of $C$ over $A$. The complex $M$ is finite and consists of
free modules, so $\hh{M\otimes_AC}\cong \hh{M\otimes_AV}$, and hence $\HH
1{M\otimes_AV}=0$, by Lemma \ref{depth sensitivity}. Now Lemma \ref{tor vanishing},
invoked with $U=M$ implies $\HH 1{{\HH 0M}\otimes_AV}=0$, so $\Tor 1A{\HH 0M}C=0$. Since
$k$ is a direct summand of $\HH 0M$, this implies $\Tor 1AkC=0$, equivalently,
$\Tor 1ACk=0$.
  
By Proposition \ref{cm:criteria}, the sequence $\bsx$ is $A$-regular, therefore $\depth
A\geq \dim A$ and $A$ is (big) Cohen-Macaulay.  Consequently, Lemma \ref{depth
  sensitivity}, now applied with $C=A$, implies $\HH iM=0$ for $i\ge 1$, as claimed. In
particular, $M$ is a finite free resolution of $\HH 0M$, so the projective dimension of
$\HH 0M$ is finite. Therefore, the projective dimension of $k$ is finite as well, since it
is a direct summand of $\HH 0M$.  Thus, $A$ is regular; see Serre \cite[Ch. IV, cor. 2,
th. 9]{Se}.  
\end{proof}

\begin{remark}
\label{compare}  
As stated in the introduction, the proof of \cite[(4.3)]{BM} is incorrect; in the
paragraph below we adopt the notation of \emph{loc.~cit.}. The problem with it is the
claim in display (1), on \cite[pg. 639]{BM}, which reads:
\[
\Tor pANC = 0 \quad \text{for all}\quad p > 0 \tag{$\ast$}
\]
This cannot hold unless we assume \emph{a priori} that the local ring $A$ is Cohen-Macaulay.

Indeed, suppose the displayed claim is true. Consider the standard change of rings spectral
sequence sitting in the first quadrant:
\[
\operatorname{E}^2_{p,q} = 
\Tor pA{\Tor qRkA}C \Longrightarrow \Tor {p+q}RkC.
\]
The edge homomorphisms in the spectral sequence give rise to the exact sequence
\[
0 = \Tor 2ANC \lra \Tor 1RkA\otimes_A C \lra \Tor 1RkC=0
\]
where the $0$ on the left holds by ($\ast$) and that on the right holds because $C$ is free over
$R$.  Thus, the middle term is $0$, which implies $\Tor 1RkA = 0$. Therefore $A$ is free
as an $R$-module, because $A$ is finitely generated over $R$, and hence $A$ is
Cohen-Macaulay.
\end{remark}

\end{document}